\documentclass [11pt,a4paper]{article}

\usepackage{amsmath}
\usepackage{amsthm}
\usepackage{graphics}
\usepackage{latexsym}
\usepackage{amssymb}
\usepackage{enumerate}
\usepackage[latin2]{inputenc}
\usepackage{color}

\newtheorem  {theorem}       {Theorem}
\newtheorem  {lemma}         {Lemma}

\newtheorem* {theorem*}      {Theorem}
\newtheorem* {lemma*}        {Lemma}
\newtheorem* {corollary*}    {Corollary}
\newtheorem* {proposition*}  {Proposition}
\newtheorem* {definition*}   {Definition}
\newtheorem* {remark*}       {Remark}
\newtheorem* {remarks*}      {Remarks}

\def \N {\mathbb N}
\def \Z {\mathbb Z}

\def \R {\mathbb R}

\def \ind{1\!\!1}

\def \dsum {\displaystyle\sum}

\def \lam   {\lambda}

\def \om {\omega}

\def \Gam   {\Gamma}

\def\wh{\widehat}

\def\beqs{\begin{eqnarray*}}
\def\eeqs{\end{eqnarray*}}
\def\beq{\begin{eqnarray}}
\def\eeq{\end{eqnarray}}
\def\beas{\begin{eqnarray*}}
\def\eeas{\end{eqnarray*}}
\def\bea{\begin{eqnarray}}
\def\eea{\end{eqnarray}}

\def \prob        {\ensuremath{\mathbf{P}}}
\def \expect      {\ensuremath{\mathbf{E}}}

\def \dt            {\, \textup{d} t}

\newcommand{\abs}[1]{\left|{#1}\right|}

\def \nnn {\nodes}                  
\def \TR  {\trees}                  
\def \ccc  {\tau}                   
\def \www   {w}                
\def \tree   {\Upsilon}             
\def \malt   {\lam^*}               
\def \lap   {\wh\rho}                  
\def \dstd   {p_{\www}}             
\def \dstt   {\pi_{\www}}           
\def \deg    {\textup{deg}}         
\def \csucs   {\eta}                

\newcommand{\und}[1] {\underline{#1}}
\newcommand{\fakt}[2] {(#1)_{#2}}

\def \repr {\xi}                    
\def \birth {\sigma}                

\def \b   {\beta}
\def \deriv {\kappa}

\newtheorem*{thmA}{Theorem A}
\newtheorem*{thmB}{Theorem B}
\newtheorem*{thmC}{Theorem C}

\def \nodes {\mathcal{N}}
\def \trees {\mathcal{G}}
\def \d {\, \textup{d}}
\def \distd   {\dstd}             
\def \distt   {\dstt}           

\newcommand{\shift}[2]   {{#1}_{\downarrow {#2}}}
\newcommand{\subtree}[2] {{#1}_{\downarrow {#2}}}
\newcommand{\nord}[1]  {\mathcal{S}({#1})}
\newcommand{\mtrees}[1] {\trees^{(#1)}}

\newcommand{\treegen}[2] {#1_{[#2]}}
\title
{Random Trees and General Branching Processes}

\author {
{\sc  Anna Rudas$^*$, B\'alint T\'oth$^*$, Benedek
Valk\'o$^{*,\circ}$} 
\\  \\
$^*$Institute of Mathematics, TU Budapest\\
$^{\circ}$Alfr\'ed R\'enyi Institute of Mathematics, Budapest
}

\begin{document}

\setlength{\baselineskip}{1.23\baselineskip}

\maketitle
\begin{abstract}

We consider a tree that grows randomly in time. Each time a new
vertex appears, it chooses exactly one of the existing vertices
and attaches to it. The probability that the new vertex chooses
vertex $x$ is proportional to $\www(\deg(x))$, a weight function of
the actual degree of $x$. The weight function
$\www:\N\to \R _+$ is the parameter of the model.

In \cite{bolriospetus} and \cite{mori1} the authors derive
the asymptotic degree distribution for a model that is
equivalent to the special case, when the weight function
is linear. The proof therein strongly relies on the linear choice
of $\www$.

Using well established results from the theory of general branching 
processes we give the asymptotical degree distribution for a wide range 
of weight functions. Moreover, we provide the asymptotic distribution
of the tree itself as seen from a randomly selected vertex. The
latter approach gives greater insight to the limiting
structure of the tree.

Our proof is robust and we believe that the method may be used to
answer several other questions related to the model.
It relies on the fact that considering the evolution of the random tree
in \emph{continuous time}, the process may be viewed as a
general branching process, this way classical results can be applied.

\end{abstract}


\section{Introduction}
\label{section:intro}

In models of randomly growing networks, the concept of
preferential attachment means that when a new vertex appears, it
is more likely that it attaches to an already existing vertex if
the latter already has more neighbors. One of the realizations of
this idea is the Barab\'asi - Albert graph. In this model, a tree
grows in discrete time-steps: at each step we add a new vertex and
connect it to an already existing vertex. This choice is made
randomly, using probabilities exactly proportional to the degree
of the existing vertices. This model reproduces phenomena observed
in certain real-world networks (see \cite{barabasialbert}), the
power-law decay of the degree sequence, for example. This was
proved in a mathematically precise way in \cite{bolriospetus} and,
independently, in \cite{mori1}.

The probability with which a newly appearing vertex chooses its
neighbor could depend on the degrees of the existing vertices in a
more general way. We define a weight function, $\www :\N \to
\R_+$, and let the probability be proportional to this function of
the degree.  In \cite{bolriospetus} and \cite{mori1}, $\www$ is
linear, their techniques strongly depend on martingales that are
apparent in the system only in this special case.

General weight functions are considered in the model of Krapivsky
and Redner \cite{krapivskyPRE}. There $\www (k) \sim k^{\gamma}$,
and non-rigorous results are obtained, showing the
different behavior for $\gamma>1$ and $\gamma\le 1$. In the first
region the limiting object does not have a non-trivial
degree sequence: a single dominant vertex appears which
is linked to almost every other vertex, the others
having only finite degree. This statement is made precise and
proved rigorously in \cite{oliveira}. (See also \cite{chung} for a related model.)
The weight functions we consider will be such that this
does not happen. 

Our main results are the following. We determine the asymptotic
distribution of the degree sequence, which equivalently gives the
limiting distribution of the degree of a (uniformly) randomly selected
vertex. We also look deeper into the structure of the tree:
we give the asymptotic distribution of the subtree under a
randomly selected vertex. Moreover, we present the asymptotic
distribution of the whole tree, seen from a randomly selected
vertex. For a general approach for asymptotic distribution of random 
subtrees of random trees, see \cite{aldous}. 
These results give greater insight to the limiting structure of the random 
tree.

The key of our method is to place the process into continuous time
in such a way that it fits in to the well-developed theory of
general branching processes (see Section \ref{section:branching}).
The asymptotic behavior of the
continuous time model gives that of the discrete time model,
as pointed out in Section \ref{section:model}. A special case of
our model (when $\www$ is linear) is equivalent to those studied
in \cite{bolriospetus} and \cite{mori1}. The robustness of our
method is apparent in the fact that it gives a.s.~results for a
wide class of weight functions, it generalizes the results mentioned
above. 

The class of weight functions for which our theorems hold are not
fully identified by the condition we make: it is sufficient, but
we do not claim that it is also necessary. In the present paper we do not intend to
analyse all the ways we believe this new approach can be used. We are
planning to apply this technique to answer many other questions related
to the randomly growing tree model.

The paper is organized as follows. We introduce the terminology and
notation in the next section, then in Section \ref{section:model} we
define the model. We state our results in Section
\ref{section:results} and show the simplifications that arise in the
linear case. After that, in Section \ref{section:branching} we give a
brief introduction to the theory of general branching processes, and
quote the theorems that we rely on. We present our proofs in Section
\ref{section:proofs}. We comment on the asymptotic growth of the tree
in the last section.


\bigskip

\section{Rooted ordered trees: terminology and notation}
\label{section:notations}

Throughout the paper it will be convenient to use genealogical
phrasing. We will consider our tree evolving as a genealogical
tree: where the individuals are the vertices and the parent child
relations are the edges of the tree. It will be also convenient to
follow the `birth orders' among children of the same parent. For
this purpose we will consider our trees always as \emph{rooted
ordered trees} (sometimes called family trees or rooted planar
trees). In the following paragraphs we introduce the (commonly
known) terminology for these trees and also some additional
notations needed for our theorems.

\subsection{Vertices, individuals, trees}
\label{subs:nodes}

We will label the vertices of a rooted ordered tree using a subset
of
\begin{equation*}
\label{def:nodes}
\nodes:=\bigcup_{n=0}^\infty\Z_+^n,
\qquad\text{ where }
\quad
\Z_+:=\{1,2,\dots\},
\quad
\Z_+^0:=\{\emptyset\}.
\end{equation*}
$\emptyset$ denotes the root of the tree, its children are
labelled with $\{1,2,\dots\}$, and in general the children of
$x=(x_1,x_2,\dots,x_k)\in\nodes$ are labelled by
$(x_1,x_2,\dots,x_k, 1), (x_1,x_2,\dots,x_k,2), \dots$. Thus if a
vertex has the label $x=(x_1,x_2,\dots,x_k)\in\nodes$ then this
means that it is  the $x_k^{\mathrm{th}}$ child of its parent,
which is the $x_{k-1}^{\mathrm{th}}$ child of its own parent and
so on. If $x=(x_1,x_2,\dots,x_k)$ and $y=(y_1,y_2,\dots,y_l)$ we
will use the shorthanded notation $x y$ for the concatenation
$(x_1,x_2,\dots,x_k,y_1,y_2,\dots,y_l)$, and with a slight abuse
of notation for $n\in \Z_+$ we use $xn$ for
$(x_1,x_2,\dots,x_k,n)$.

We will identify a rooted ordered tree with the set of labels of
its vertices, since this already contains the necessary
information about the edges. It is clear that a $G\subset \nnn$
may represent a rooted ordered tree if and only if $\emptyset \in
G$ and for each $(x_1,x_2,\dots,x_k)\in G$ we have
$(x_1,x_2,\dots,x_{k-1})\in G$ as well as
$(x_1,x_2,\dots,x_k-1)\in G$, if $x_k>1$.

The set of finite rooted ordered trees will be denoted by
$\trees$. We think about $G\in\trees$ as an oriented tree with edges
pointing from parents to children. The \emph{degree} of a vertex
$x\in G$ is just the number of its children in $G$:
\begin{equation*}
\label{def:degree} \deg(x,G):=\max{\{n\in \Z_+:xn\in G\}}.
\end{equation*}
The \emph{$n^{\textup{th}}$ generation} of $G\in\trees$ is
\begin{equation*}
\label{def:nthgen}
\treegen{G}{n}:=\{x\in G: |x|=n\}, \quad n\ge0,
\end{equation*}
where $|x|=n$ iff $x\in \Z_+^{n}$.

The \emph{$n^{\textup{th}}$ ancestor of $x=(x_1,x_2,\dots,x_k)\in
\nnn$} with $k\ge n$ is $x^n=(x_1,x_2,\dots,x_{k-n})$ if $k>n$
and $x^n=\emptyset$ if $k=n$.

The \emph{subtree} rooted at a vertex $x\in G$ is:
\begin{equation*}
\label{def:subtree} \subtree{G}{x}:= \{y: xy\in G\},
\end{equation*}
this is just the progeny of $x$ viewed as a rooted ordered tree.
Also, (again with a slight abuse of notations) for an
$x=(x_1,x_2,\dots,x_n)\in \nnn$ with $\abs{x}=n\ge k$ we use the
notation $\subtree{x}{k}=(x_{n-k+1},x_{n-k+2},\dots,x_n)$. This
would be the new label given to $x\in G$ in the subtree
$\subtree{G}{x^k}$.

Consider a $G\in \trees$. An ordering
$s=(s_0,s_1,\dots,s_{\abs{G}-1})$ of the elements of $G$ is called
\emph{historical} if it gives a possible 'birth order' of the vertices
in $G$,
formally if for each $0\le i \le {\abs{G}-1}$ we have
$\{s_0,s_1,\dots,s_i\}\in \trees$. The set of all historical orderings
of $G\in\trees$ will be denoted $\nord{G}$.
For a fixed  $s\in \nord{G}$ the rooted ordered trees
\begin{equation}
\label{def:history} G(s,i):=\{s_0,s_1,\dots,s_i\}\subset G
\end{equation}
give the evolution of $G$ in this historical ordering $s$.

\subsection{Random trees}

Throughout the paper we will use Greek letters to denote random
elements (of various distributions) selected from $\nodes$ and
$\trees$:
\begin{equation*}
\label{not:rndnodesandtrees}
\xi, \zeta, \dots \in{\nodes}, \qquad
\Gamma,\Upsilon, \dots\in{\trees}
\end{equation*}

Our results will deal with some asymptotic properties of a
randomly chosen vertex in a certain random tree. We will investigate
the asymptotic  distribution of its degree, its progeny and also
the progeny of its the $k^{\textup{th}}$ ancestor. In order to
study the latter object, we introduce \emph{rooted ordered trees
with a marked vertex in generation $k$}:
\begin{equation*}
\label{def:markedtrees}
\mtrees{k}:=\{(G,u)\in\trees\times\Z^k\,:\,u\in \treegen{G}{k}\}.
\end{equation*}
$\mtrees{0}$ is identified with $\trees$, since generation
0 consists of only the root, $\emptyset$. We can use the elements of
$\mtrees{k}$ to describe the progeny of the $k^{\textup{th}}$
ancestor of a random vertex: $G$ is an ordered tree rooted in the
$k^{\textup{th}}$ ancestor of the selected point and
$u\in \treegen{G}{k}$ is the
position of the random vertex in this tree. Clearly, if $(G,u)$
describes the progeny of the $k^{\textup{th}}$ ancestor, then for
$0\le l \le k$ the progeny of the $l^{\textup{th}}$ ancestor is
described by $(\subtree{G}{u^{l}},\subtree{u}{k-l})$.

Thus if $\pi^{(k)}$ is a distribution on $\mtrees{k}$ which
describes the progeny of the $k^{\textup{th}}$ ancestor of a
chosen vertex, then, if $l<k$, the distribution of the progeny of the
$l^{\textup{th}}$ ancestor is:
\[
\pi^{(k,l)}(H,v):=\pi^{(k)} \left( \left\{ (G,u)\in \mtrees{k}:
\subtree{G}{u^{l}}=G, v=\subtree{u}{k-l} \right\}\right)
\]
The sequence $\pi^{(k)}$ of probability measures on $\mtrees{k}$,
$k=0,1,2,\dots$ is called \emph{consistent} if for any $0\le l\le
k$, the identity $\pi^{(l)} = \pi^{(k,l)}$ holds.

Without presenting the precise formulation, it is clear that a
consistent sequence $\pi^{(k)}$ of probability measures on
$\mtrees{k}$ gives full insight to the limiting structure of the
tree as seen from a random vertex, see Remark after Theorem
\ref{thm:main}.

We call a probability measure $\pi$ on $\trees$  \emph{steady} if
\begin{equation}
\label{def:steadypm} \sum_{H\in\trees} \pi(H)\sum_{x\in
\treegen{H}{1}} \ind\{\subtree{H}{x}=G\}= \pi(G).
\end{equation}
It is easy to check that in this case, for any $k=1,2,\dots$, the
similar identity
\begin{equation*}
\label{steady2} \sum_{H\in\trees} \pi(H)\sum_{x\in \treegen{H}{k}}
\ind\{\subtree{H}{x}=G\}= \pi(G)
\end{equation*}
follows.
Equivalently, for any bounded function $\varphi:{\trees}\to\R$ and any
$k=0,1,2,\dots$,
\begin{equation*}
\label{steady3} \expect \big(\abs{\treegen{\Gamma}{k}}
\varphi(\subtree{\Gamma}{\xi})\big) =
\expect \big(\varphi(\Gamma)\big),
\end{equation*}
where $\Gamma$ is a random element of $\trees$ with distribution
$\prob(\Gamma=G)=\pi(G)$, and on the left hand side $\xi$ is a random
vertex selected uniformly from the $k$-th generation of $\Gamma$. (We
don't have to worry about the fact that $\treegen{\Gamma}{k}$ may be
empty, since in that case the expression $\abs{\treegen{\Gamma}{k}}
\varphi(\subtree{\Gamma}{\xi})$ is automatically 0.) Immediate
consequences of this property are that the expected size of the
$k^{\text{th}}$ generation is 1 for any $k\in \N$ (choose $\varphi$
identically 1), and therefore the expected size of the whole tree is
infinite.

\emph{Backward extensions:} Given a steady probability measure
$\pi$ on $\trees$, define the probability measures $\pi^{(k)}$ on
$\mtrees{k}$, $k=0,1,2\dots$,  by
\begin{equation*}
\label{steadyextension} \pi^{(k)}(G,u):=\pi(G)
\end{equation*}
One can easily check that, due to the steadiness of the
distribution $\pi$, the sequence of probability measures
$\pi^{(k)}$ on $\mtrees{k}$, $k=0,1,2\dots$, is consistent.


\section{The Random Tree Model}
\label{section:model}

\subsection{Discrete Time Model}

Given a  weight function $\www:\N\to\R_+$, let us define the following 
discrete time Markov chain $\tree_d$ on the countable state space $\trees$, 
with initial state $\tree_d (0)= \{\emptyset\}$. If for $n\ge 0$ we have 
$\tree_d (n)=G$, then for a vertex $x \in G$ let $k:=\deg(x,G)+1$.
Using this notation, let the transition probabilities be 
\begin{equation*}
\prob(\tree_d (n+1)= G\cup \{xk\})=\frac{\www(\deg (x,G))} {\sum_{y\in
G}\www(\deg(y,G))}.
\end{equation*}

In other words, at each time step a new vertex appears, and attaches to
exactly one already existing vertex. If the tree at the appropriate
time is $G$, then the probability of choosing vertex $x$
in the tree $G$ is proportional to $\www(\deg (x,G))$.

\subsection{Continuous Time Model}
Given a  weight function $\www:\N\to\R_+$,
let $X(t)$ be a Markovian pure birth process with $X(0)=0$ and birth
rates
\begin{equation*}
\label{birthrates}
\prob\big(X(t+\dt)=k+1\,\big|\,X(t)=k\big)=\www(k)\dt +o(\dt).
\end{equation*}
Let $\rho:[0,\infty)\mapsto(0,\infty]$ be the density of the point
process corresponding to the pure birth process $X(t)$ and
$\wh\rho:(0,\infty)\to(0,\infty]$ the (formal) Laplace transform of
$\rho$:
\begin{equation}
\label{def:laprho}
\wh\rho(\lam):=
\int_0^\infty e^{-\lambda t}\rho(t)\dt=
\dsum_{n=0}^{\infty} \prod_{i=0}^{n-1}
\frac{\www(i)}{\lam+\www(i)}.
\end{equation}
The rightmost expression of $\wh\rho(\lambda)$ is easily computed,
given the fact that the intervals between successive jumps of
$X(t)$ are independent exponentially distributed random variables
of parameters $\www(0), \www(1), \www(2), \dots$
respectively. Let
\begin{equation}
\label{def:undlam}
\underline{\lambda}:=\inf\{\lambda>0: \wh\rho(\lambda)<\infty\}.
\end{equation}
Throughout this paper we impose the following condition on the weight
function $\www$:
\begin{equation*}
\label{condM}
\tag{M}
\lim_{\lambda\searrow\underline\lambda}\wh\rho(\lambda)>1.
\end{equation*}

We are now ready to define our randomly growing  tree $\tree(t)$
which will be a continuous time, time-homogeneous Markov chain on
the countable state space $\trees$, with initial state
$\tree(0)=\{\emptyset\}$.

The jump rates are the following: if for a $t\ge 0$ we have
$\tree(t)=G$ then the process may jump to $G\cup \{xk\}$ with rate
$\www(\deg(x,G))$ where $x\in G$ and $k=\deg(x,G)+1$. This means
that each existing vertex $x\in \tree(t)$ `gives birth to a child'
with rate $\www (\deg (x, \tree (t)))$ independently of the
others.

Note that condition \eqref{condM} implies
\begin{equation*}
\label{noexpl}
\sum_{k=0}^\infty \frac1{\www(k)}=\infty
\end{equation*}
and hence it follows that the Markov chain $\tree(t)$ is well defined
for $t\in[0,\infty)$, it does not blow up in finite time. A rigorous proof of this statement follows from the connection with general branching processes (see Section \ref{section:branching}) for which the related statement is derived in \cite{jagers}.

We define the \emph{total weight} of a tree $G\in\trees$ as
\begin{equation*}
\label{def:osszsuly}
W(G):=\sum_{x\in G} \www(\deg(x,G)).
\end{equation*}
Described in other words, the Markov chain $\tree(t)$ evolves as
follows:  assuming $\tree(t-)=G$,
at time $t$ a new vertex is added to it  with total rate
$W(G)$ which is attached with an oriented edge (pointing
towards the newly added vertex) to the already existing vertex $x\in G$
with probability
\begin{equation*}
\label{def:kapcs_prob} \frac{\www(\deg (x,G))} {\sum_{y\in
G}\www(\deg(y,G))}.
\end{equation*}
Therefore, if we only look at our process at the stopping times
when a new vertex is just added to the randomly growing tree:
\begin{equation}
\label{def:T_n} T_n:=\inf \{t: \abs{\tree (t)} =n+1\}
\end{equation}
then we get the discrete time model: $\tree (T_n)$ has the same
distribution as $\tree_d (n)$, the discrete time model at time $n$.

It will sometimes be convenient  to refer to the vertices in the
order of their birth, not their genealogical code: let
$\{\csucs_k\}=\tree(T_k)\setminus\tree(T_k-)$ denote the vertex
that appeared at $T_k$. Of course we will always have $\csucs
_0=\emptyset$ and $\csucs_1=1$.


\section{Results}
\label{section:results}

\subsection{Statement of results}
\label{subs:results}

From  condition \eqref{condM}  it follows that the equation
\begin{equation*}
\label{rho=1}
\wh\rho(\lambda)=1
\end{equation*}
has a unique root $\lambda^*$.

Now we are ready to state our first theorem.

\begin{theorem}
\label{thm:meta}
Consider a weight function $\www$ satisfying condition \eqref{condM}
and let $\malt$ be defined as above. Consider a bounded function
$\varphi:\trees\to\R$. Then the following limit holds almost surely:
\begin{equation*}
\label{eq:thm_meta}
\lim_{t\to \infty}\frac{1}{\abs{\tree(t)}}
\sum_{x\in \tree(t)} \varphi(\subtree{\tree(t)}{x})=
\malt
\int_0^{\infty} e^{-\malt\, t} \expect\big(\varphi(\tree(t))\big)
\dt.
\end{equation*}
\end{theorem}

From Theorem \ref{thm:meta} several statements follow, regarding
the asymptotic behavior of our random tree  \emph{as seen from a
randomly selected vertex} $\zeta$, chosen uniformly from
$\tree(t)$. As typical examples we  determine the asymptotic
distribution of the number of children, respectively, that of the
whole subtree under the randomly chosen vertex, its
$k^{\textup{th}}$ ancestor, respectively. That is: the asymptotic
distribution of $\deg(\zeta,\tree(t))\in\N$, 
$\subtree{\tree(t)}{\zeta}\in\trees$ and
$(\subtree{\tree(t)}{\zeta^{(k)}},\subtree{\zeta}{k})
\in\mtrees{k}$.

In order to formulate these consequences of Theorem \ref{thm:meta}
we need to introduce some more notation. Let $G\in\trees$ and one
of its historical orderings
$s=(s_0,s_1,\dots,s_{\abs{G}-1})\in\nord{G}$
be fixed. The historical sequence of total weights are defined as%
\beq \label{def:histW} W(G,s,i):=W(G(s,i)) \eeq%
for $0 \le i \le \abs{G}-1$ while the respective weights of the
appearing vertices are defined as
\begin{eqnarray}
\label{def:histw}
w(G,s,i):=w\left(\deg\left((s_i)^1,G(s,i-1)\right)\right).
\end{eqnarray}
for $1 \le i \le \abs{G}-1$. Since
$\deg\left((s_i)^1,G(s,i-1)\right)$ is the degree of $s_i$'s
parent just before $s_i$ appeared,  $w(G,s,i)$ is the rate with
which our random tree process jumps from $G(s,i-1)$ to $G(s,i)$.

Given the weight function $\www:\N\to\R_+$ satisfying condition
(\ref{condM}) and $\malt$ defined as before define
\begin{eqnarray}
p_{\www}(k)&:=& \frac{\malt}{\malt+\www(k)} \prod_{i=0}^{k-1}
\frac{\www(i)}{\malt+\www(i)},\label{def:p_om}
\\[5pt]
\pi_{\www}(G)&:=& \sum_{s\in\nord{G}}\label{def:pi_om}
\frac{\malt}{\malt+W(G)} \prod_{i=0}^{\abs{G}-2}
\frac{w(G,s,i+1)}{\malt+W(G,s,i)}.
\end{eqnarray}

\begin{theorem}
\label{thm:main}
Consider a weight function $\www$ which satisfies condition (M) and
let $\malt$ be defined as before. Then the following limits hold
almost surely:

\begin{enumerate}[(a)]

\item
For any fixed $k\in \N$
\begin{equation*}
\label{eq:thm_degree}
\hskip-3mm
\lim_{t\to\infty}
\frac{\abs{\{x\in \tree(t): \deg(x,\tree(t))=k\}}}{\abs{\tree(t)}}
=
\distd(k).
\end{equation*}

\item
For any fixed $G\in\trees$
\begin{equation*}
\label{eq:thm_graph}
\lim_{t\to\infty}
\frac{
\abs{\{x\in\tree(t):\subtree{\tree(t)}{x}=G\}}}{\abs{\tree(t)}}
=
\distt(G).
\end{equation*}

\item
For any fixed $(G,u)\in\mtrees{k}$
\begin{equation*}
\label{eq:thm_familytree} \lim_{t\to\infty} \frac{
\abs{\{x\in\tree(t):
(\subtree{\tree(t)}{x^{(k)}},\shift{x}{k})=(G,u)\}}}{\abs{\tree(t)}}
= \distt(G).
\end{equation*}

\end{enumerate}

Furthermore, the functions $\distd, \distt$ are probability
distributions on $\N$ and $\trees$, respectively, and $\distt$ is
steady (i.e.~identity \eqref{def:steadypm} holds).

\end{theorem}

{\bf \noindent Remarks.}\\
\noindent \textbf{1.} Parts (a), (b) and (c) of Theorem
\ref{thm:main}, in turn, give more and more information about the
asymptotic shape of the randomly growing tree $\tree(t)$, as seen
from a random vertex $\zeta$ chosen with uniform distribution. Part
(a) identifies the a.s.~limit as $t\rightarrow\infty$, of the
degree distribution of $\zeta$. Part (b) identifies the a.s.~limit
as $t\rightarrow\infty$, of the distribution of the progeny of
 $\zeta$. Finally, part (c) does the same for the distribution of
the progeny of the $k^{\textup{th}}$ ancestor of the randomly
selected vertex with the position of this vertex marked.

\noindent \textbf{2.} From part (c) it is easy to derive the
asymptotic distribution of the progeny of the $k^{\textup{th}}$
ancestor of the randomly selected vertex (as a rooted ordered tree
without any marked vertices):
\[
\lim_{t\to\infty} \frac{ \abs{\{x\in\tree(t):
\subtree{\tree(t)}{x^{(k)}}=G\}}}{\abs{\tree(t)}} =\distt(G)
\abs{\treegen{G}{k}}.
\]
The limit is the size-biased version of $\distt(G)$, with the
biasing done by the size of the $k^{\textup{th}}$ generation.

\noindent \textbf{3.} Since the distribution $\distt$ is steady, part
(c) identifies the asymptotic distribution of the \emph{whole family
tree} of the randomly selected vertex $\zeta$ (relatives of arbitrary
degree included). Hence asymptotically, as $t\rightarrow \infty$, the
tree $\tree(t)$ viewed from a random vertex $\zeta$ will have the
following structure (we omit the precise formulation): \\
-- there exists an infinite path of ancestors $\zeta^1,\zeta^2,\zeta^3,
\dots$ `going back in time',\\
-- we have finite ordered random trees rooted at each vertex of this path,\\
-- the tree rooted at $\zeta^k$ with the position of $\zeta$ marked on it
has distribution $\distt^{(k)}$ on $\mtrees{k}$ where
$\distt^{(k)}(G,u)=\distt(G)$.

\subsection{Linear weight function}
\label{subs:linear}

In the linear case $\www(k)=\alpha k+\beta$ $(\alpha, \beta >0)$ all computations
are rather
explicit. In this case the asymptotic degree distribution $\distd$
(computed in \cite{bolriospetus}, \cite{mori1}) is reproduced, of
course. But even in this simplest case the asymptotic distribution
$\distt$ of the subtree under a randomly selected vertex seems to be
new.

For sake of completeness, in the rest of this section we perform
these (explicit and straightforward) computations for the linear
case. Multiplying the rate function with a positive constant only
means the rescaling of time in our model thus it is enough to
consider $\www(k)=k+\beta$ (with $\beta>0$). In this case it is
straightforward to compute that condition (M) holds,
$\lap(\lam)=\frac{\beta}{\lam-1}$, $\und{\lam}=1$ and
$\malt=1+\beta$. Thus both Theorems \ref{thm:meta} and
\ref{thm:main} hold.

For the asymptotic degree distribution we get
\begin{equation*}
\label{lindegdistr}
\distd(k)=
(1+\b) \frac{\fakt{k-1+\b}{k}}{\fakt{k+1+2\b}{k+1}},
\end{equation*}
where we used the shorthanded notation
\begin{equation*}
\label{def:fakt}
\fakt{x}{k}:=\prod_{i=0}^{k-1} (x-i)=\frac{\Gam(x+1)}{\Gam(x-k+1)},
\qquad
k=0,1,2,\dots.
\end{equation*}

For the calculation of $\distt(G)$ first we show that the sum
which defines it contains identical elements. In order to avoid
heavy notation, during the following computations we will use
$n:=\abs{G}-1$ and $\deg(x)$ instead of $\deg(x,G)$.

Clearly, for any $s\in \nord{G}$
\begin{equation*}
\label{linwprod} \prod_{i=0}^{n-1} w(G,s,i+1) = \prod_{x\in G}
\left( \prod_{j=0}^{\deg(x)-1} \www(j) \right) =\prod_{x\in G}
\fakt{\deg(x)-1+\b}{\deg(x)}.
\end{equation*}
(Actually, the first equality holds for every weight function
$\www$.)
It is also easy to see that for any $G\in \trees$
\begin{equation*}
\label{linW} W(G)= \sum_{x\in G} (\deg(x)+\b)=\abs{G}(1+\b)-1,
\end{equation*}
thus for any $s\in \nord{G}$
\begin{equation*}
\label{lincompute1} \frac{\malt}{\malt+W(G)} \prod_{i=0}^{n-1}
\frac{1}{\malt+W(G,s,i)}=
\frac{1}{(1+\b)^n {\fakt{n+2-(1+\b)^{-1}}{n+1}}}.
\end{equation*}
Therefore
\begin{equation*}
\label{lingraphdistr} \distt(G)= \abs{\nord{G}} \frac{\prod_{x\in
G} \fakt{\deg(x)-1+\b}{\deg(x)}}
     {{(1+\b)^n}{\fakt{n+2-(1+\b)^{-1}}{n+1}}}.
\end{equation*}

In the $\b=1$ case (i.e. if we consider random tree proposed in
\cite{barabasialbert}) the previous calculations give
\begin{equation*}
\label{beta1lindegdistr}
\distd(k)=\frac{4}{(k+1)(k+2)(k+3)}
\end{equation*}
and
\begin{equation*}
\label{beta1lingraphdistr} \distt(G)=
\frac{2\abs{\nord{G}}}{(2\abs{G}+1)!!} \prod_{x\in G} \deg(x)!.
\end{equation*}
\bigskip

The value of $\abs{\nord{G}}$ cannot be written as the function of
degrees of $G$ only, but one can compute it using the values
$\abs{\subtree{G}{x}}$ for $x\in G$. For a given $G$ and
$x=(x_1,x_2,\dots,x_n)\in G$ let us introduce the following
notations (these will not be used in the other parts of the
paper): \beqs &B(x):=\{y\in G:
y=(x_1,x_2,\dots,x_{n-1},k), k>x_n \},& \\
&a(x):=\max(\abs{\subtree{G}{x}}-1,1),\qquad
b(x):=\max\left(\sum_{y\in B(x)}\abs{\subtree{G}{x}},1 \right).&
\eeqs%
It is a simple exercise to prove that for a $G\in \trees$ with
$\abs{G}>1$
\[
\abs{\nord{G}}=(\abs{G}-2)!{\prod_{ x\in G, x\neq \emptyset}
a(x)^{-1} b(x)^{-1}}.
\]
The proof is left to the reader.


\section{Branching Processes}
\label{section:branching}

The random tree model, defined in continuous time, has the big
advantage that it fits into the framework of the well-established
theory of general branching processes. We give a brief
introduction to the fundamentals and state the theorems that we
rely on in our proofs. We do not give a broad survey on the most
general types of branching processes here, we choose to focus on
the results which may be applied to our process. For more details
see the monograph \cite{jagers} or the papers \cite{jagersnerman},
\cite{nerman}, \cite{olofsson} and the references therein. For a
survey on branching processes, trees and superprocesses, see
\cite{legall}.

In the case of a general branching process, there is a population
in which each individual reproduces at ages according to
i.i.d.~copies of a random point process $\repr$ on $[0,\infty)$.
We denote by $\repr (t)$ the $\repr$-measure of $[0,t]$, this the
random number of children an individual has up to time $t$.


The individuals in the population are labelled with the elements
of $\nnn$, as described in Section \ref{section:model} (see
\eqref{def:nodes}). The basic probability space is
\[
(\Omega, \mathcal{A}, P)= \prod _{x\in \nnn}(\Omega_x,
\mathcal{A}_x, P_x),
\]
where $(\Omega_x, \mathcal{A}_x, P_x)$ are identical spaces on
which $\repr _x$ are distributed like $\repr$.

For each $x\in \nnn$ there is a $.\shift{}{x}$ shift defined on $\Omega$
by
\[
(\omega_{\shift{}{x}})_y=\omega_{xy},
\]
in plain words, $\omega_{\shift{}{x}}$ is the life of the progeny of
$x$, regarding $x$ as the ancestor.

The birth times
$\birth _x$ of the individuals are defined in the obvious way:
$\birth _0 =0$ and if $x'=xn$ with $n\in \Z_+$ then %
\beq \label{def:sigma}\birth _{x'}=\birth_x + \inf \{t: \repr_x
(t)\ge n\}. \eeq

The branching process is often counted by a random characteristic,
this can be any real-valued process $\{\Phi:\R\times \Omega \to \R\}$.
For each individual $x$, $\Phi_x$ is defined by
\[
\Phi_x(t,\omega)=\Phi(t, \omega_{\shift{}{x}}),
\]
in plain words $\Phi _x(t)$ denotes
the value of $\Phi$ evaluated on the progeny of $x$, regarding $x$
as the ancestor, at the time when $x$ is of age $t$. We can think
about $\Phi _x(t)$ as a `score' given to $x$ when its age is $t$.
With this,
\[
Z_t^{\Phi} = \sum _{x\in \nnn} \Phi _x (t-\birth _x)
\]
is the branching process counted by the random characteristic
$\Phi$ (the `total score' of the population at time $t$).

\emph{For our applications we only consider random characteristics
which are 0 for $t<0$ and equal to a bounded
\underline{deterministic} function of the rooted tree for $t\ge 0$. }

This means that only those individuals contribute to $Z_t^{\Phi}$
which are born up to time $t$ and their contribution is a
deterministic function of their progeny tree. (Random
characteristics may be defined in a more general way, see
e.g.~\cite{jagers}, \cite{jagersnerman}.) One of the important
examples is $\Phi(t) =\ind \{t\ge 0\}$ when $Z_t^{\Phi}$
is just the total number of individuals 
born up to time $t$.

The Laplace-transform of $\textup{d} \repr(t)$ will be of great
importance,
we denote this random variable by:%
\beq%
\wh{\repr}(\lam):=\int_0^{\infty} e^{-\lam t} \d
\repr(t).\label{def:dxi}
\eeq%
We shall be interested in \emph{supercritical, Malthusian}
processes, meaning that there exists a finite $0<\malt<\infty$
(the so-called Malthusian parameter) for which
\beq\label{def:malthus}%
\expect\, \wh{\repr}(\malt)=1,
\eeq%
and also %
\beq\label{def:deriv} \deriv=-\partial_{\lambda} \left(\expect\,
\wh{\repr}(\lambda)\right) \Big|_{\lambda=\malt}=\expect
\int_0^{\infty} t e^{-\malt t} \d \repr(t)<\infty. \eeq (The last
property means that the process is Malthusian and the first means
that it is supercritical.)
Also, we require the reproduction to be non-lattice, which means
that the jumps of $\repr(t)$ cannot be supported by any lattice
$\{0, d, 2d, \dots\}$, $d>0$ with probability one.

We quote here a weaker form of Theorem 6.3 from \cite{nerman},
using its extension which appears in Section 7 of the same paper.
This way the conditions of the original theorem are fulfilled
automatically.

\begin{thmA}[Nerman, \cite{nerman}]\label{thm:"6.3"}
Consider a supercritical, Malthusian branching process with
Malthusian parameter $\malt$, counted by two random
characteristics $\Phi(t)$ and $\Psi(t)$ which have the properties
described above (i.e.~they are 0 for $t<0$ and a deterministic
bounded
function of the progeny tree for $t\ge 0$). Suppose that 
there exists a
$\und{\lambda}< \malt$ for which
\[\expect \,\wh{\repr}(\und{\lam})<\infty.\]
Then 
almost surely \beq \frac{Z_t^{\Phi}}{Z_t^{\Psi}} \to \frac{\wh{
\Phi}(\malt)}{\wh{\Psi}(\malt)} \quad \text{as }t\to \infty, \eeq
where $\wh{\Phi}(\lambda) = \int\limits_0^{\infty}\exp (-\lambda
t) \expect (\Phi (t))\dt$.
\end{thmA}

For a reader interested in how the Malthusian parameter and
$\wh{\Phi}(\lambda)$ play a role in the theory, we give a short
indication. This part may be skipped without any confusion.

The key observation is the so-called \textit{basic decomposition},
namely that \beq \label{eq:basic_decomp}Z_t^{\Phi}=\Phi
(t)+\sum_{j \in \N} Z_{t-\birth_j}^{\Phi_j}. \eeq The advantage of
this formula is that if we know the sequence $(\birth_j)_{j\in
\Z_{+}}$ (all the birth-times of the children of the root), then
$Z_{t-\birth_j}^{\Phi_j}$ has the same conditional distribution as
$Z_{t-\birth_j}^{\Phi}$.

Therefore, using the notation $m_t^{\Phi}=\expect Z_t^{\Phi}$,
taking expectation on both sides of \eqref{eq:basic_decomp} in two
steps (first conditionally on $(\birth_j)_{j \in \Z_{+}}$, then
taking expectation regarding $(\birth_j)_{j \in \Z_{+}})$, we get

\beq m_t^{\Phi}=\expect (\Phi (t))+ \int_0^t m_{t-s}^{\Phi} \d
\mu(s), \eeq where we used the notation $\mu(t)=\expect \repr(t)$.
 Taking the
Laplace transform of both sides, \beq \wh{m}(\lambda)=\wh
{\Phi}(\lambda)+\wh{m}(\lambda) \wh{\mu}(\lambda), \eeq so
formally \beq
\wh{m}(\lambda)=\frac{\wh{\Phi}(\lambda)}{1-\wh{\mu}(\lambda)}.
\eeq Note that $\wh{\mu}(\lambda)=\expect \,\wh{\repr}(\lam)$.
This shows that if there is a positive interval below $\malt$
where the Laplace transform is finite, then
$1/(1-\wh{\mu}(\lambda))$ has a simple pole at $\malt$ (it is easy
to check that $\wh{\mu}'(\malt)<0$ and $\wh{\mu}''(\malt)>0$). So
taking series expansion and inverse Laplace transform results that
\beq m_t^{\Phi}=\frac1{\deriv}{\wh{\Phi}(\malt)} \, e^{\malt
t}+o(e^{\malt t}), \eeq where $\deriv$ is the finite constant
defined in (\ref{def:deriv}). This means that the ratio of the
expectations of $Z_t^{\Phi}$ and $Z_t^{\Psi}$ indeed tends to
$\frac{\wh{\Phi}(\malt)}{\wh{\Psi}(\malt)}$. To get the almost
sure convergence of $\frac{Z_t^{\Phi}}{Z_t^{\Psi}}$ to the same
limit needs a lot more work and of course its proof is more
elaborate (see \cite{nerman} for example).


\section{Proofs}
\label{section:proofs}

\begin{proof}[Proof of Theorem \ref{thm:meta}]
Consider the continuous time branching process where the
reproduction process $\repr(t)$ is the Markovian  pure birth
process $X(t)$, with rate function $\www$, described at the
beginning of Section \ref{section:model}. 

Clearly, the time-evolution of the population has the same
distribution as the evolution of the continuous time Random Tree
Model corresponding to the weight function $\www$. The vertices
are the respective individuals and edges are the
parent-child relations.

It is also not hard to see that the function $\expect
\,\wh{\repr}(\lam)$ for the branching process is  the same as
$\wh\rho(\lam)$ which means that by condition (M) we may apply
Theorem A with appropriate random characteristics.
Given any bounded function $\varphi:\trees\to\R$, setting the
characteristics $\Phi, \Psi$  as
$\Phi(t):=\varphi(\tree(t))\,  \ind\{t\ge 0\}$ and
$\Psi(t):=\ind\{t\ge 0\}$ we get exactly the statement of Theorem
\ref{thm:meta}.
\end{proof}

\begin{proof}[Proof of Theorem \ref{thm:main}]
(a)
Apply Theorem \ref{thm:meta} with the function
\begin{equation*}
\label{def:Thetadeg}
\varphi(G):=
\ind\{\deg(\emptyset,G)=k\}.
\end{equation*}
This gives that
\begin{equation*}
\label{eq:proof_thm2a}
\lim_{t\to\infty}
\frac{\abs{\{x\in \tree(t): \deg(x,\tree(t))=k\}}}{\abs{\tree(t)}}
=
\malt \int_0^{\infty} e^{-\malt\, t} \,
\prob\big(\deg(\emptyset, \tree(t))=k\big) \dt,
\end{equation*}
almost surely.
By the definition of $\birth_k$ (see (\ref{def:sigma})):%
\begin{equation*}
\label{degcompute1}
\prob\big(\deg(\emptyset,\tree(t))=k\big)
=
\prob\big(\birth_k<t\big)-\prob\big(\birth_{k+1}<t\big).
\end{equation*}
Since
\begin{equation*}
\label{degcompute2}
\malt \int_0^{\infty} e^{-\malt\, t}\,  \prob\big(\birth_k<t\big)
\dt=\expect \big(e^{-\malt \birth_k}\big),
\end{equation*}
and $\birth_k$ is the sum of independent exponentially distributed
random variables with parameters $\www(0), \www(1), \dots,
\www(k-1)$, we get
\begin{equation*}
\label{degcompute3}
\malt \int_0^{\infty} e^{-\malt\, t}
\prob\big(\deg(\emptyset,\tree(t))=k\big)\dt
=
\frac{\malt}{\malt+\www(k)} \prod_{i=0}^{k-1}
\frac{\www(i)}{\malt+\www(i)}.
\end{equation*}
This completes the proof of part (a) of the Theorem.

Using the identity
\begin{equation*}
\distd(k)=\malt \int_0^{\infty} e^{-\malt\, t} \,
\prob\big(\deg(\emptyset,\tree(t))=k\big) \dt,
\end{equation*}
and the fact that $\abs{\tree(t)}$ is finite for every
$t$ with probability 1
 it is straightforward to prove that $\distd$ is
indeed a probability distribution on $\N$.

\medskip

\noindent
(b)
Let $G\in \TR$ be fixed and denote
$n=\abs{G}-1$. We apply Theorem \ref{thm:meta} with
$\varphi(H)=\ind( H=G )$.  We need to compute
\begin{equation*}
\label{eq:proof_thm2b_1}
\malt \int_0^{\infty} e^{-\malt\, t}\,
\prob\big(\tree(t)=G\big) \dt.
\end{equation*}
Consider the following random stopping times:
\begin{eqnarray*}
\label{def:stop}
\tau_{G}
&:=&
\sup\{t\ge0:\tree(t)\subseteq G\},
\\
\notag
\tau'_{G}
&:=&
\sup\{t\ge0:\tree(t)\subsetneq G \}.
\end{eqnarray*}
That is: $\tau_G$ is the birth time of the first vertex
not in $G$, while $\tau'_G$ is the minimum of
$\tau_G$ and the time when we first have $\tree(t)=G$, if
the latter ever happens. Since
\begin{equation*}
\label{graphcompute1}
\prob\big(\tree(t)=G\big)
=
\prob\big(\tree(t)\subseteq G\big)-\prob\big(\tree(t)\subsetneq G\big)
=
\prob\big(t<\tau_G\big)-\prob\big(t<\tau'_G\big),
\end{equation*}
we get that
\begin{eqnarray*}
\label{graphcompute2}
\malt \int_0^{\infty} e^{-\malt\, t}\,
\prob\big(\tree(t)=G\big) \dt
&=&
\expect\big(e^{-\malt \tau'_{G}}-e^{-\malt \tau_{G}}\big)
\\
\notag
&=&
\expect \big((e^{-\malt \tau'_{G}}-e^{-\malt\tau_{G}})
\ind\{\tau'_{G}<\tau_{G}\}\big).
\end{eqnarray*}
Note  that by the definition we always have $\tau'_{G}\le\tau_{G}$.
The event $\{\tau'_{G}<\tau_{G}\}$ means that there is a $t$
when $\tree(t)= G$. On this event $\tau'_{G}$ gives the time
when we first have $\tree(t)= G$ and $\tau_{G}$ gives the
appearance of the next vertex. Given the event
$\{\tau'_{G}<\tau_{G}\}$, the conditional distribution of
$\tau_{G}-\tau'_{G}$ is exponential with parameter
$W(G)$
and it is (conditionally) independent of $\tau'_{G}$.
This leads to
\begin{equation*}
\label{graphcompute3}
\expect \big((e^{-\malt\tau'_{G}}-e^{-\malt\tau_{G}})
\ind\{\tau'_{G}<\tau_{G}\}\big)
=
\frac{\malt}{\malt+W(G)} \,
\expect\big(e^{-\malt \tau'_{G}}\ind\{\tau'_{G}<\tau_{G}\}\big).
\end{equation*}
Now, it is clear that the following two events are actually the same
\begin{equation*}
\label{graphcompute4}
\big\{\tau'_{G}<\tau_{G}\big\}\,=\,
\big\{(\csucs_0,\dots,\csucs_{n}) =
   (s_0,\dots, s_n) \textup{ for some }
   s\in \nord{G}\big\}.
\end{equation*}
This implies that
\begin{eqnarray*}
\label{graphcompute5}
&&
\expect \big(e^{-\malt \tau'_{G}}\ind\{\tau'_{G}<\tau_{G}\}\big)
\\
\notag
&&
\hskip20mm
=
\sum_{s\in\nord{G}}
\expect \big(e^{-\malt T_n}
\ind\{(\csucs_0,\dots,\csucs_n)=(s_0,\dots,s_n)\} \big).
\end{eqnarray*}
For $T_n$ see \eqref{def:T_n} and the definition below it.
Given $s\in \nord{G}$ fixed
\begin{equation*}
\label{graphcompute6}
\prob\big((\csucs_0,\dots,\csucs_n)=(s_0,\dots,s_n)\big)=
\prod_{i=0}^{n-1}\frac{w(G,s,i+1)}{W(G,s,i)}.
\end{equation*}
(See \eqref{def:kapcs_prob}, \eqref{def:histW} and
\eqref{def:histw} for the definitions.) Also, if $s\in
\nord{G}$ is fixed then \emph{conditionally on the
  event}
$\{(\csucs_0,\dots,\csucs_n)=(s_0,\dots,s_n)\}$
the random variables $T_{k+1}-T_{k}$, $k=0,1, \dots,n-1$,
are independent and exponentially distributed with parameters
$W(G,s,k)$, $k=0,1,\dots,n-1$, respectively.
This is an easy exercise: it may be
proved by using the `lack of memory' of the exponential
distribution and the fact that the minimum of independent
exponentially distributed random variables with parameters
$\nu_1,\nu_2,\dots,\nu_l$ is also exponentially distributed with
parameter $\sum_{i=1}^l \nu_i$. Hence  it is straightforward to get
\begin{equation*}
\label{graphcompute7}
\expect \big(e^{-\malt T_n}
\ind\{(\csucs_0,\dots,\csucs_n)=(s_0,\dots,s_n)\}\big)
=
\prod_{i=0}^{n-1}\frac{w(G,s,i+1)}{\malt+W(G,s,i)}.
\end{equation*}
Collecting our previous calculations part (b) of Theorem
\ref{thm:main} follows.

Using similar considerations as in the end of the proof of part
(a) it is apparent that $\distt$ is a probability
distribution on $\trees$.
\medskip

\noindent (c) This is straightforward since for any $H\in \trees$ and
$(G,u) \in \mtrees{k}$ we have%
\beqs%
\abs{\{x\in H: (\subtree{H}{x^{(k)}},\shift{x}{k})=(G,u)\}}=
{\abs{\{x\in H:\subtree{H}{x}=G\}}}. \eeqs%
The statement now follows from part (b).
\medskip

The only thing left to prove is that $\distt$ satisfies
(\ref{def:steadypm}), i.e.~it is steady. First observe, that if
$G_0\in \trees$ is fixed and $\zeta$ is a uniformly chosen random
vertex in $G_0$ then the distribution of
$\Gamma:=\subtree{(G_0)}{\zeta}$ (which is a probability
distribution on $\trees$) is steady. (This follows by simple
counting.) Equation (\ref{def:steadypm}) is linear in $\pi$,
therefore mixtures of steady distributions are also steady. Thus,
if $\zeta$ is a uniformly chosen random vertex in $\tree(t)$ then
the distribution of $\subtree{\tree(t)}{\zeta}$ (which is a
\emph{random} probability distribution on $\trees$) is also
steady. By part (b) with probability one these distributions
converge (in distribution) to $\distt$ and from this an easy
consideration shows that $\distt$ must satisfy
(\ref{def:steadypm}).
\end{proof}


\section{Asymptotic growth}
\label{section:growth}

In our theorems we determined the asymptotic ratio of vertices in
$\tree(t)$ satisfying certain properties. It is also natural to
ask if one can prove results about the asymptotic number of the
respective vertices. As we have seen, this essentially requires to
study the asymptotic behavior  of $Z_t^{\Phi}$ for a suitable
random characteristic $\Phi$. This has been done in the framework
of general branching processes, we shall give a short overview of
the relevant results.

We have already seen that \beq \expect (e^{-\malt t}Z_t^{\Phi})\to
\frac1{\deriv} {\hat{\Phi}(\malt)}, \eeq where $\kappa$ is the
constant defined in (\ref{def:deriv}). Thus we need to divide
$Z_t^{\Phi}$ by $e^{\malt t}$ to get something non-trivial. Let us
quote a weaker form of Theorem 5.4 of \cite{nerman}.

\begin{thmB}[Nerman, \cite{nerman}]\label{thm:"5.4"}

Consider a supercritical, Malthusian branching process with
Malthusian parameter $\malt$. Suppose that condition (M) holds and
$\Phi$ is a random characteristic with properties described in
Section \ref{section:results}. Then almost surely
\[
e^{-\malt t}Z_t^{\Phi}\to \frac1{\deriv}
{\hat{\Phi}(\malt)}\Theta, \quad \text{as }t \to \infty,
\]
where $\Theta$ is a random variable not depending on $\Phi$.
\end{thmB}

The necessary and sufficient condition for the random variable
$\Theta$ to be a.s.~positive is the so-called $x \log x$ property
of the reproduction process $\repr$:
\begin{equation*}
\label{condL} \tag{L} \expect (\wh{\repr}(\malt)\log ^+
\wh{\repr}(\malt))<\infty.
\end{equation*}
We quote Theorem 5.3 of \cite{jagersnerman}.

\begin{thmC}[Jagers-Nerman, \cite{jagersnerman}]
Consider a supercritical, Malthusian branching process with
Malthusian parameter $\malt$. If condition (\ref{condL}) holds
then $\Theta>0$ a.s.~and $\expect (\Theta)=1$; otherwise
$\Theta=0$ a.s.
\end{thmC}

{\bf \noindent Remark.} This theorem is the generalization of the
Kesten-Stigum theorem, which states this fact for Galton-Watson
processes (see \cite{kestenstigum}).

Theorem B applies for the random tree  model if
condition (M) is fulfilled. We do not intend  to identify the
necessary and sufficient condition on the weight function $\www$
which would guarantee that the corresponding reproduction process
possesses property (L). Still it is worth pointing out that if
$w(k)\to \infty$ as $t\to \infty$, then this property holds, this by
Theorem C, $\Theta$ is a.s.~positive.

\begin{lemma}
If a weight function $\www$ satisfies condition (M) and
$\www(n)\rightarrow \infty$, as $n\to \infty$,
then the corresponding branching process satisfies condition (L).
\end{lemma}

\begin{proof}
We will prove the existence of the second moment of
$\wh{\repr}(\malt)$ from which condition (L) trivially follows.
Since $\wh{\repr}(\malt)=\sum_{i=1}^{\infty} e^{-\malt \birth_k}$,
we need %
\beq \expect \left(\sum_{k=1}^{\infty} e^{-\malt
\birth_k}\right)^2<\infty. \label{eq:2nd_mom}\eeq%
The random variables  $\birth_{k+1}-\birth_{k}$ are independent
exponentials for $k=0,1,2,\dots$ with parameters $\www(0),
\www(1), \dots$, respectively, thus a simple computation yields
that the expression in (\ref{eq:2nd_mom}) is equal to
\beqs%
-\lap(2\malt)+2 \sum_{i=0}^{\infty} \sum_{j=0}^{i}\left(
\prod_{l=0}^j \frac{\www(l)}{2 \malt +\www(l)} \prod_{l=j+1}^{i}
\frac{\www(l)}{\malt +\www(l)}\right).
\eeqs%
Transforming the double sum on the right, we get%
\beqs%
\sum_{i=0}^{\infty} \left(\prod_{l=0}^i \frac{\www(l)}{\malt
+\www(l)} \right) \sum_{j=0}^{i} \prod_{l=0}^j \frac{\malt
+\www(l)}{2 \malt +\www(l)} \le  \sup_{ i} \sum_{j=0}^{i}
\prod_{l=0}^j \frac{\malt +\www(l)}{2 \malt +\www(l)}
\eeqs%
where we also used $\lap(\malt)=1$. On the other hand, %
\beqs%
\sum_{j=0}^{i} \prod_{l=0}^j \frac{\malt +\www(l)}{2 \malt
+\www(l)}=\sum_{j=0}^{i} \prod_{l=0}^j \frac{\www(l)}{\und{\lam}
+\www(l)}\prod_{l=0}^j \frac{(\und{\lam} +\www(l))(\malt
+\www(l))}{\www(l)(2 \malt +\www(l))}.
\eeqs%
Since $\www(l)\rightarrow\infty$ and $\und{\lam}<\malt$, we have
\beqs \frac{(\und{\lam} +\www(l))(\malt +\www(l))}{\www(l)(2 \malt
+\www(l))}<1 \eeqs if $l$ is large enough.
This leads to %
\[ \sum_{j=0}^{i} \prod_{l=0}^j \frac{\malt +\www(l)}{2 \malt
+\www(l)} <K \, \sum_{j=0}^{i} \prod_{l=0}^j
\frac{\www(l)}{\und{\lam} +\www(l)} < K \lap(\und{\lam})<\infty,\]
by condition (M) which completes the proof of the lemma.
\end{proof}

\def \tot {N}

The distribution of $\Theta$ is usually hard to determine from the
weight function $\www$, however, one can characterize its moment
generating function $\varphi(u)=\expect e^{-u \Theta}$. Using the
idea of the basic decomposition (\ref{eq:basic_decomp}) one can
write the following equation for $f(u,t):=\expect e^{-u
\abs{\tree(t)}}$
\[
f(u,t)= e^{-u}\expect \prod _{j\ge 1} f(u, t-\ccc _j).
\]
By Theorem 4  \beq\varphi(u)=\lim_{t\to \infty}f(ue^{-\malt
t},t)\label{eq:lim}\eeq
which gives%
\beq\label{eq:momgen} \varphi(u)= \expect  \prod _{j\ge 1}
\varphi(ue^{-\malt \ccc
_j}). \eeq%
It can be proved that this equation characterizes $\varphi$ as
there is no other bounded function satisfying it with a right
derivative $-1$ at 0. (See Theorem 6.8.3 in \cite{jagers}.)

In the linear case ($\www (k)=\alpha k+\beta$) the distribution of $\Theta$
may be explicitly calculated. Using the fact that $\abs{\tree(t)}$ is
now a Markov process (which is an easy exercise to prove) and applying
standard martingale techniques one can derive a partial differential
equation for $f(u,t)$. Solving this pde and using the limit
(\ref{eq:lim}) we get the momentum generating function of $\om$ from
which one can identify the distribution of $\Theta$ as a Gamma
distribution with parameters $(\frac{\b}{\b +\alpha},\frac{\b}{\b +\alpha})$.


\noindent{\bf\large Acknowledgements:}
\\
We are thankful to J-F.~Le Gall for having called our attention
to his review paper \cite{legall} and the reference
\cite{jagersnerman} therein.
This work was partially supported by
the Hungarian Scientific Research Fund (OTKA) grants no.~T037685, TS40719 
and K60708.


\vfill

\hbox{\sc
\vbox{\noindent
\hsize66mm
Anna Rudas \& B\'alint T\'oth\\
Institute of Mathematics\\
Technical University Budapest\\
Egry J\'ozsef u. 1.\\
H-1111 Budapest, Hungary\\
{\tt \{rudasa,balint\}{@}math.bme.hu}
}
\hskip20mm
\vbox{\noindent
\hsize53mm
Benedek Valk\'o\\
R\'enyi Inst. of Math.\\
Hung. Acad. of Sciences\\
Re\'altanoda u. 13-15\\
H-1053 Budapest, Hungary\\
{\tt valko{@}renyi.hu}
}
}

\end{document}